\documentclass[12pt]{article}
\usepackage{graphicx}
\usepackage{array}
\usepackage{fullpage}
\usepackage[small, hang, bf]{caption}
\begin{document}

\title{Calculation of Exact Estimators by Integration Over the Surface
  of an 
  $n$-Dimensional Sphere}

\author{
Anthony J. Webster\thanks{email: dr.anthony.webster@gmail.com}\\
%anthony.webster@ccfe.ac.uk}
%\normalsize{Abingdon, Oxfordshire, United Kingdom}
%\normalsize{United Kingdom Atomic Energy Authority,}\\
%\normalsize{Culham Science Centre, Abingdon, Oxon, OX14 3DB.}
}

\maketitle

\begin{abstract}

This paper reconsiders the problem of calculating the expected 
set of probabilities $\langle p_i \rangle$, given the observed set of
items $\{ m_i \}$, that are distributed among $n$ bins with an
(unknown) set of probabilities $\{ p_i\}$ for  being placed in the
$i$th bin.  
The problem is often 
formulated using Bayes theorem and the
multinomial distribution, along with a constant prior for the values
of the $p_i$, leading to a Dirichlet distribution for the $\{ p_i
\}$. 
The moments of the $p_i$ can then be calculated exactly. 
Here a new approach is suggested for the calculation of the moments,
that uses a change of variables that 
reduces the problem to an integration over 
a portion of 
the surface of an
$n$-dimensional sphere.  
This greatly simplifies the calculation by allowing a straightforward
integration over $(n-1)$ independent variables, with the constraints
on the set of $p_i$ being automatically satisfied.  
For the Dirichlet and similar distributions the problem simplifies even
further, with the resulting integrals subsequently factorising,
allowing their easy evaluation in terms of Beta functions.  
A proof by induction confirms existing calculations for the moments. 
The advantage of the approach presented here is that the methods and
results apply with minimum or no modifications to numerical
calculations that involve more complicated 
distributions or non-constant prior distributions, for which cases the
numerical calculations will be greatly simplified.   

\end{abstract}

\section{Introduction}

Many problems involve placing $N$ objects into $n$ bins, with
probabilities $p_i$ for the object being placed into the $i$th bin.  
Given the values of the set of $\{ p_i \}$, the probability density 
$P(m_1,m_2, ... , m_n|p_1,p_2, ... , p_n)$ for the distribution of the
set of $\{ m_i\}$ objects can be calculated, and is well-know as the
multinomial distribution, 
\begin{equation}
P(m_1,m_2, ... , m_n | p_1, p_2, ... , p_n) = \frac{N!}{m_1!m_2! ... m_n!} 
\Pi_{i=1}^n p_i^{m_i}
\end{equation} 
with the constraint that $\sum_{i=1}^n p_i = 1$ 
and $\sum_{i=1}^n m_i =N$.  
Bayes theorem, $P(A|B)P(B)=P(B|A)P(A)$ requires, 
\begin{equation}
P(p_1,p_2, ... , p_n | m_1, m_2, ... ,m_n)P(m_1, m_2, ... , m_n) = 
P(m_1,m_2, ... , m_n | p_1, p_2, ... , p_n)P(p_1, p_2, ... , p_n) 
\end{equation} 
that in principle allows us to calculate 
$P(p_1,p_2, ... , p_n | m_1, m_2, ... m_n)$, the probability of the
set of probabilities $\{ p_i \}$ with $i=1$ to $i=n$, given the
observed set of $\{ m_i \}$.  
Often in such problems, $P(p_1, p_2, ... , p_n)$ is taken to be
constant, and $P(m_1,m_2, ... , m_n)$ is chosen to ensure that
$P(p_1,p_2, ... , p_n)$ is correctly normalised \cite{Jaynes}.   
%%%%%%%%
Applying this approach to the multinomial distribution, leads to a
Dirichlet distribution, for which exactly calculated moments can be
obtained. 
%%%(e.g. see page 181 of \cite{CB}).
%%%%
A recent approach to this problem by Friedman \cite{F}, relied on
an identity discovered by Gauss that involves the integral representation
of the hypergeometric distribution. 
The same is true of a recent exact calculation that
corrects conjectured  but widely used mark-recapture
estimates \cite{Webster}, this 
and the coincidental timing of its revision on arXiv are 
what brought this problem to the author's attention. 

%.   cacluatio is what brought this problem to
%the author's attention. 
%Interestingly, 
%his conjecture relies on an identity discovered by
%Gauss that involves the integral representation of the hypergeometric
%distribution, and the same is true for a recent exact calculation that
%corrects a conjectured  but widely used mark-recapture
%estimate\cite{Webster}.  
%% used an alternative method showed how this
%approach can be used to calculate $\langle p_i \rangle$ and $\langle
%p_i^2 \rangle - \langle p_i \rangle^2$ exactly for $n=2$ to $n=5$, and
%observed a pattern in his calculation that led him to conjecture a
%more general result.  
%Interestingly, his conjecture relies on an identity discovered by
%Gauss that involves the integral representation of the hypergeometric
%distribution, and the same is true for a recent exact calculation that
%corrects a conjectured  but widely used mark-recapture
%estimate\cite{Webster}.    
%It was the coincidental appearance of a revision to that
%paper\cite{Webster} on arXiv at the same time as the appearance of
%Friedman's paper that led the author to consider this problem.  

Here an alternative method of calculation is considered. 
%%%the problem studied by Friedman\cite{F} is reconsidered. 
I suggest a change of variables that elegantly leads to a simple
calculation for the moments of the $\{ p_i \}$, 
%%$\langle p_i^a \rangle$, 
and confirms  existing results.  
%%conjectured by Friedman\cite{F}.   
%The results confirm those conjectured by Friedman\cite{F}. 
The advantage of the method is that it can be applied very generally,
and allows comparatively 
straightforward numerical integrations for the most general situations
when analytical solutions may not be possible.   
%%%%
The crux of the problem is the integration of 
a function 
%%over all a set of $\{ p_i \}$, 
%%the product 
%%$\Pi_{i=1}^n p_i^{m_i}$ 
over all possible values of $p_i$ between $0$
and $1$, subject to the constraint of $\sum_{i=1}^n p_i =1$. 
This appears in many situations, the specific case considered here is
the product $\Pi_{i=1}^n p_i^{m_i}$ that arises in the Binomial,
Multinomial, and Dirichlet distributions for example. 

%for example the ? distribution has, 
%\begin{equation}
%asdfasdfds
%\end{equation} 
%and the ? distribtuion has,
%\begin{equation}
%asdfadsf
%\end{equation}
%%%both have this form, with. 
%The solution to this generic problem is described next.

\section{The Calculation}

Consider the integration of 
%%%%%%%We would like to integrate 
the product $\Pi_{i=1}^n p_i^{m_i}$, over
all sets of 
values of the $p_i$, subject to the constraints of $0\leq p_i \leq 1$ for
all $i$, and $\sum_{i=1}^n p_i =1$.  
%%%
In Casella and Berger \cite{CB}, the moments are obtained by a
delightful 
trick (page 181), that simplifies the problem to integration over a
binomial distribution.  
In Friedman \cite{F} the integral is accomplished by a nested set of
integrals, each of 
which depends on the calculation of the integrals within it, with for
$n=3$ for example,  
\begin{equation}\label{IF}
I_3 = \int_{p_1=0}^1 dp_1 \int_{p_2=0}^{1-p_1} dp_2 p_1^{m_1}
p_2^{m_2} \left( 1 - p_1 - p_2 \right)^{m_3}  
\end{equation}
where $\sum_{i=1}^3 p_i = 1$ has been used to write $p_3=1-p_1-p_2$. 
%%%%%
Here I start in a similar way, writing, 
\begin{equation}\label{e1}
\Pi_{i=1}^n p_i^{m_i} = 
\left( 1 - \sum_{i=1}^{n-1} p_i \right)^{m_n} 
\Pi_{i=1}^{n-1} p_i^{m_i} 
\end{equation}
that for $n=3$ is $p_1^{m_1}p_2^{m_2}(1-p_1-p_2)^{m_3}$. 
%%%%%%%%%%%
%Here we start in a similar way, writing, 
%\begin{equation}\label{e1}
%\Pi_{i=1}^n p_i^{m_i} = \left(1 - \sum_{i=1}^{n-1} p_i \right)^{m_n} 
%\Pi_{i=1}^{n-1} p_i^{m_i} 
%\end{equation}
%that for $n=3$ is $p_1^{m_1}p_2^{m_2}(1-p_1-p_2)^{m_3}$. 
%%%%%%%%%%
Eq. \ref{e1} recognises that the constraint of $\sum_{i=1}^n p_i=1$
leads to $(n-1)$ free parameters, or $2$ free parameters for $n=3$.  
For a radius of $r=1$ the $n$-dimensional polar co-ordinates are: 
\begin{equation}\label{ndpolar}
\begin{array}{rl}
x_1(n) &= \cos \theta_1
\\
x_2(n) &= \sin \theta_1 \cos \theta_2
\\
x_3(n) &= \sin \theta_1 \sin \theta_2 \cos \theta_3
\\
... & ...
\\
x_{n-1}(n) &= \sin \theta_1 \sin \theta_2 ... \sin \theta_{n-2} \cos
\theta_{n-1}  
\\ 
x_{n}(n) &= \sin \theta_1 \sin \theta_2 ... \sin \theta_{n-2} \sin
\theta_{n-1}  
\end{array}
\end{equation}
Notice that $x_i(n)$ and $x_i(n)^2$ will vary continuously between $0$
and $1$ as the set of $\theta_i$ are varied continuously between $0$
and $\pi/2$.  
Also notice that $\sum_{i=1}^n x_i(n)^2 =1$, and consequently that 
$x_n(n)^2 = 1 - \sum_{i=1}^{n-1} x_i(n)^2$. 
Therefore the substitutions of $p_1=x_1(n)^2$, $p_2=x_2(n)^2$, ... ,
$p_{n-1}=x_{n-1}(n)^2$, will ensure that $\sum_{i=1}^n p_i =1$, and
integrals over $\theta_i$ from $\theta_i=0$ to $\pi/2$ will allow
$p_i$ to vary continuously over all values between $0$ and $1$.  

Note that the constraint of $\sum_{i=1}^n p_i =1$ leads to $(n-1)$
free parameters, that after the change of variables, correspond to the
set of $\theta_i$ with $i=1$ to $(n-1)$. Also note that although we
are using polar co-ordinates in $n$ dimensions, because we have set
$r=1$, there are only $(n-1)$ free parameters. 

The Jacobian of the co-ordinate transformation is $J= \left| \partial
x_i(n)^2 / \partial \theta_j \right|$.  
%%% 
Notice from Eq. \ref{ndpolar} that $\partial x_i(n)^2/\partial
\theta_j=0$ for $j>i$.   
Consequently the determinant has zeros above the diagonal, and will
evaluate easily to give $J= \Pi_{i=1}^{n-1} 
\left| \partial x_i(n)^2/\partial \theta_i \right|$.  

Before proceeding to the general case, consider again the case with
$n=3$, for which case, 
\begin{equation}
\begin{array}{l}
x_1(3) = \cos \theta_1
\\
x_2(3) = \sin \theta_1 \cos \theta_2
\\
x_3(3) = \sin \theta_1 \sin \theta_2
\end{array}
\end{equation}
The product $\left( 1 - \sum_{i=1}^{n-1} p_i \right)^{m_n}
\Pi_{i=1}^{n-1} p_i^{m_i}$ becomes, after the change of variables,  
\begin{equation}\label{Intg1}
\begin{array}{ll}
\left( 1 - p_1 - p_2 \right)^{m_3} p_1^{m_1} p_2^{m_2} &= 
\left( \sin^2 \theta_1 \sin^2 \theta_2 \right)^{m_3} 
\left( \cos^2 \theta_1 \right)^{m_1}
\left( \sin^2 \theta_1 \cos^2 \theta_2 \right)^{m_2}
\\
&= 
\left( \cos^{2m_1}\theta_1 \sin^{2(m_2+m_3)} \theta_1 \right) 
\left( \cos^{2m_2}\theta_2 \sin^{2m_3} \theta_2 \right)  
\end{array}
\end{equation}
The Jacobian is, 
%\begin{equation}
%%\begin{array}{l}
%J = \left|
%\begin{array}{ll}
%\frac{\partial x_1(3)^2}{\partial \theta_1}
%& \frac{\partial x_1(3)^2}{\partial \theta_2}
%\\
%\frac{\partial x_2(3)^2}{\partial \theta_1}
%& \frac{\partial x_2(3)^2}{\partial \theta_2
%\end{array}
%\right|
%\end{equation}
\begin{equation}\label{J1}
\begin{array}{ll}
J
&=
\left|
\begin{array}{ll}
-2\cos \theta_1 \sin \theta_1 & 0
\\
2 \sin \theta_1 \cos \theta_1 \cos^2 \theta_2 & 
-2 \sin^2 \theta_1 \sin \theta_2 \cos \theta_2 
\end{array}
\right|
%\end{equation}
%%%%%%%%%%%%%%%%%%%%
%\begin{equation}
\\
&= \left( 2 \cos \theta_1 \sin^3 \theta_1 \right)
\left( 2 \sin \theta_2 \cos \theta_2 \right) 
\end{array}
\end{equation}
Therefore using Eqs. \ref{Intg1} and \ref{J1} the integral in Eq
\ref{IF} can be equivalently calculated from,  
\begin{equation}
I_3 = \int_0^{\pi/2} d\theta_1 \int_0^{\pi/2} d\theta_2 
\left( \cos^{2m_1}\theta_1 \sin^{2(m_2+m_3)} \theta_1 \right) 
\left( \cos^{2m_2}\theta_2 \sin^{2m_3} \theta_2 \right)  
\left( 2 \cos \theta_1 \sin^3 \theta_1 \right)
\left( 2 \sin \theta_2 \cos \theta_2 \right)
\end{equation} 
This integral factorises into, 
\begin{equation}\label{I3ind}
I_3 = \left(
2 \int_0^{\pi/2} d\theta_1 \cos^{2(m_1+1)-1}\theta_1 
\sin^{2(m_2+m_3+2)-1} \theta_1 
\right)
\left(
2 \int_0^{\pi/2} d\theta_2
\cos^{2(m_2+1)-1} \theta_2 \sin^{2(m_3+1)-1} \theta_2 
\right)
\end{equation}
the above Eq. \ref{I3ind} will be used as a starting point for a proof
by induction for the general case later.  

Many readers will immediately recognise the integrals as Beta
functions, and it is well known that,  
\begin{equation}\label{eql12}
2 \int_0^{\pi/2} d\theta \cos^{2m-1} \theta \sin^{2n-1} \theta =
\mbox{B}(m,n) =  
\frac{\Gamma(m)\Gamma(n)}{\Gamma(m+n)}
\end{equation} 
Consequently $I_3$ is easily evaluated as, 
\begin{equation}
I_3 = 
\frac{\Gamma(m_1+1)\Gamma(m_2+m_3+2)}{\Gamma(m_1+m_2+m_3+3)}
\frac{\Gamma(m_2+1)\Gamma(m_3+1)}{\Gamma(m_2+m_3+2)}
\end{equation}
Cancelling terms and writing in terms of factorials this gives, 
\begin{equation}
I_3 = \frac{m_1!m_2!m_3!}{(m_1+m_2+m_3+2)!} 
\end{equation}
For non-integer $m_i$ Eq. \ref{eql12} must be left written in terms of
Gamma functions.

If we now wish to calculate $\langle p_1 \rangle$ for example, we simply
need to evaluate
$I_3(m_1+1,m_2,m_3)/I_3(m_1,m_2,m_3)=(m_1+1)/(m_1+m_2+m_3+3)=(m_1+1)/(N+3)$
with $N=m_1+m_2+m_3$, as found by Friedman. 
Other moments are easily calculated in a similar way. 

For the general case, consider the formulae,
\begin{equation}\label{Ind0a}
I_n = \int_0^{\pi/2} d\theta_1 \int_0^{\pi/2} d\theta_2 ...
\int_0^{\pi/2} d\theta_{n-1} \Pi_{j=1}^{n-1} K_j(n) 
\end{equation}
\begin{equation}\label{Ind0b}
K_j(n) = 2 \cos^{2(m_j+1)-1} (\theta_j) 
\sin^{2\sum_{l=j+1}^n (1 + m_l) -1 } (\theta_j) 
\end{equation}
where I note that $\sum_{l=j+1}^n(1+m_l)=(n-j)+\sum_{l=j+1}^n m_l$,
and the dependency on $n$ of $K_j(n)$ is through the upper limit in
the sum.  
%%%%%
Note that Eqs. \ref{Ind0a} and \ref{Ind0b} are true for $n=3$, as can
be seen by comparison with Eq. \ref{I3ind}.  
%%%%
I will assume this is true for $n=k$ then show that this implies 
it is true for $n=k+1$, and consequently for all $k\geq 3$
by induction.  

Firstly consider the integral with $n=k$.
For $n=k$ the change of variables is, 
\begin{equation}
\begin{array}{rrl}
p_1 &= x_1(n)^2 &= \cos^2 \theta_1
\\
p_2 &= x_2(n)^2 &= \sin^2 \theta_1 \cos^2 \theta_2
\\
p_3 &= x_3(n)^2 &= \sin^2 \theta_1 \sin^2 \theta_2 \cos^2 \theta_3
\\
... & ...
\\
p_{k-1} &= x_{k-1}(n)^2 &= \sin^2 \theta_1 \sin^2 \theta_2 ... \sin^2
\theta_{k-2}  \cos^2 \theta_{k-1}  
\\
p_k &= x_{k}(n)^2 &= \sin^2 \theta_1 \sin^2 \theta_2 ... \sin^2
\theta_{k-2} \sin^2  \theta_{k-1}  
\end{array}
\end{equation}
and the integrand is $\Pi_{i=1}^k p_i^{m_i}$, with a Jacobian that as
noted previously, 
simplifies to $J=\Pi_{i=1}^{k-1} \left| \partial (x_i(k)^2)/\partial
\theta_i \right|$.  
This gives the integral $I_k$ as, 
\begin{equation}\label{Ik}
I_k = \int_0^{\pi/2} d\theta_1 ... \int_0^{\pi/2} d\theta_{k-1} 
\Pi_{i=1}^k x_i(k)^{2m_i} \Pi_{j=1}^{k-1} 
\left| \partial x_j(k)^2/\partial \theta_j \right| 
\end{equation}
Now consider $n=k+1$, for which the change of variables is, 
\begin{equation}
\begin{array}{rrl}
p_1 &= x_1(n)^2 &= \cos^2 \theta_1
\\
p_2 &= x_2(n)^2 &= \sin^2 \theta_1 \cos^2 \theta_2
\\
p_3 &= x_3(n)^2 &= \sin^2 \theta_1 \sin^2 \theta_2 \cos^2 \theta_3
\\
... & ...
\\
p_{k-1} &= x_{k-1}(n)^2 &= \sin^2 \theta_1 \sin^2 \theta_2 ... \sin^2
\theta_{k-2}  \cos^2 \theta_{k-1}  
\\ 
p_k &= x_{k}(n)^2 &= \sin^2 \theta_1 \sin^2 \theta_2 ... \sin^2
\theta_{k-2} \sin^2  \theta_{k-1} \cos^2 \theta_{k}  
\\
p_{k+1} &= x_{k+1}(n)^2 &= \sin^2 \theta_1 \sin^2 \theta_2 ... \sin^2
\theta_{k-2} \sin^2  \theta_{k-1} \sin^2 \theta_k 
\end{array}
\end{equation}
and the integral $I_{k+1}$ is, 
\begin{equation}\label{Ik+1}
I_{k+1} = \int_0^{\pi/2} d\theta_1 ... \int_0^{\pi/2} d\theta_{k} 
\Pi_{i=1}^{k+1} x_i(k+1)^{2m_i} \Pi_{j=1}^{k} 
\left| \partial x_j(k+1)^2/\partial \theta_j \right|
\end{equation}
Now notice that for $i=1$ to $i=(k-1)$, $x_i(k)=x_i(k+1)$. 
For $i=k$, $x_k(k+1)=x_k(k)\cos^2\theta_k$. 
Therefore,
\begin{equation}
\begin{array}{ll}
\Pi_{i=1}^{k+1} x_i(k+1)^{2m_i} &= \Pi_{i=1}^k x_i(k)^{2m_i} 
\cos^{2m_k} (\theta_k) x_{k+1}(k+1)
\\
&= \Pi_{i=1}^k x_i(k)^{2m_i} 
\cos^{2m_k} (\theta_k) 
\sin^{2m_{k+1}} (\theta_1) 
\sin^{2m_{k+1}} (\theta_2) ...
\sin^{2m_{k+1}} (\theta_k) 
\end{array}
\end{equation} 
Similarly the Jacobian can be written as,
\begin{equation}
\begin{array}{ll}
J &= \Pi_{i=1}^{k} 
\left| \frac{\partial}{\partial \theta_i} \left(
x_i(k+1)^2 \right)  \right| 
\\
&= \left| \frac{\partial}{\partial \theta_k}  \left( x_k(k+1)^2
\right) \right|  
\Pi_{i=1}^{k-1} 
\left| \frac{\partial}{\partial \theta_i} \left( x_i(k)^2 \right)
\right| 
\\
&= -2 
\sin^2 (\theta_1) 
\sin^2 (\theta_2) ... 
\sin^2 (\theta_{k-1}) 
\sin (\theta_k) 
\cos (\theta_k)  
\Pi_{i=1}^{k-1} 
\left| \frac{\partial}{\partial \theta_i} \left( x_i(k)^2 \right)
\right| 
\end{array}
\end{equation}
Therefore we have,
\begin{equation}\label{newIk}
\begin{array}{c}
I_{k+1} = \int_0^{\pi/2} d\theta_1 ... \int_0^{\pi/2} d\theta_{k-1} 
\int_0^{\pi/2} d\theta_k 
\Pi_{i=1}^k x_i(k)^2 
\Pi_{i=1}^{k-1} 
\left| \frac{\partial x_i(k)^2}{\partial \theta_i} \right|
\\
\sin^{2(m_{k+1}+1)} (\theta_1) ... \sin^{2(m_{k+1}+1)} (\theta_{k-1}) 
2 \cos^{2(m_{k+1}+1)-1} (\theta_k) \sin^{2(m_{k+1}+1)-1} (\theta_k) 
\end{array}
\end{equation}
Comparing Eq. \ref{Ik} with the assumption of Eq. \ref{Ind0a},  
%%and \ref{Ind0b}, 
we find,  
\begin{equation}
\Pi_{i=1}^k x_i(k)^{2m_i} \Pi_{i=1}^{k-1} 
\left| \frac{\partial x_i(k)^2}{\partial \theta_i}  \right| 
= \Pi_{i=1}^{k-1} K_j(k) 
\end{equation}
with $K_j(k)$ given by Eq. \ref{Ind0b}.
%%%%%%%%%%%%%%%%%55
%%%%
%%%
Under this assumption the integrand of Eq. \ref{newIk} can be written as, 
\begin{equation}\label{Ig23}
\left[
2 \cos^{2(m_{k+1}+1)-1} (\theta_k) \sin^{2(m_{k+1}+1)-1} (\theta_k) 
\right]
\Pi_{i=1}^{k-1} \left[
K_j(k) \sin^{2(m_{k+1}+1)} (\theta_j) 
\right]
\end{equation}
Note that, 
\begin{equation}
\begin{array}{ll}
K_j(k)\sin^{2(m_{k+1}+1)} (\theta_j) &=
2 \cos^{2(m_j+1)-1} (\theta_j) 
\sin^{2\sum_{l=j+1}^{k+1} (1+m_l) -1} (\theta_j)
\\
&= K_j(k+1) \mbox{  for  } 1\leq j \leq (k-1)
\end{array}
\end{equation}
The extra factor in Eq. \ref{Ig23} is, 
\begin{equation}
2 \cos^{2(m_{k+1}+1)-1} (\theta_k) 
\sin^{2(m_{k+1}+1)-1} (\theta_k) = K_k(k+1) 
\end{equation}
Therefore we have, 
\begin{equation}\label{Ikp1}
I_{k+1} = \int_0^{\pi/2} d\theta_1 ... \int_0^{\pi/2} d\theta_k 
\Pi_{i=1}^k K_i(k+1) 
\end{equation}
which is just Eq. \ref{Ind0a} with $n=(k+1)$, and $K_i(k+1)$ as given
by Eq. \ref{Ind0b}.  
Since we've shown Eq. \ref{Ikp1} to be true for $n=3$ and that its
truth for $n=k$ implies it to be true for $n=(k+1)$, then by induction
Eqs. \ref{Ind0a} and \ref{Ind0b} are true for all $n\geq 3$.

Eq. \ref{Ikp1} is easy to evaluate. 
Because $\theta_i$ only appears in $K_i(k+1)$, the integral factors
into, 
\begin{equation}
I_{k+1} = \Pi_{i=1}^k \int_0^{\pi/2} d\theta_i K_i(k+1) 
\end{equation}
Noting Eq. \ref{Ind0b} for $K_i(k+1)$, each of the integrals can be
recognised as a Beta function, with, 
\begin{equation}
\begin{array}{ll}
  \int_0^{\pi/2} d\theta_i K_i(k+1) &=
2 \int_0^{\pi/2} \cos^{2(m_i+1)-1} (\theta_i)
\sin^{2\sum_{l=i+1}^{k+1} (1+m_l) -1 } (\theta_i)
\\
&= \frac{\Gamma\left(m_i+1\right) \Gamma\left(\sum_{l=i+1}^{k+1}
  (1+m_l)\right)} 
{\Gamma\left(\sum_{l=j}^{k+1} (1+m_l)\right)} 
\end{array}
\end{equation}
where in the denominator of the last line we used 
$m_i+1 + \sum_{l=i+1}^{k+1} (1+m_l) = \sum_{l=i}^{k+1} (1+m_l)$. 
To obtain an explicit value for the integral, now we simply need to
multiply out the terms, with, 
\begin{equation}
\begin{array}{c}
I_{k+1} = 
\frac{\Gamma(m_1+1) \Gamma\left( \sum_{l=2}^{k+1} (1+m_l) \right)}
{\Gamma \left( \sum_{l=1}^{k+1} (1+m_l) \right)}
\times 
\frac{\Gamma(m_2+1) \Gamma\left( \sum_{l=3}^{k+1} (1+m_l) \right)}
{\Gamma \left( \sum_{l=2}^{k+1} (1+m_l) \right)} \times ...
\\
%\times ... \times
%\\
... \times  
\frac{\Gamma(m_{k-1}+1) \Gamma\left( m_k + m_{k+1} + 2 \right)}
{\Gamma \left( m_{k-1} + m_k + m_{k+1} +3 \right)}
\times 
\frac{\Gamma(m_{k}+1) \Gamma\left( m_{k+1} + 1 \right)}
{\Gamma \left( m_k + m_{k+1} + 2 \right)}
\end{array}
\end{equation}
Cancelling successive terms, leaves, 
\begin{equation}
I_{k+1} = \frac{ \Gamma(m_1+1) \Gamma(m_2+1)
  ... \Gamma(m_k+1)\Gamma(m_{k+1}+1)} 
{\Gamma\left( \sum_{l=1}^{k+1} (1+m_l) \right)}
\end{equation}
which when written in terms of factorials and $N = \sum_{l=1}^{k+1}
m_l$, gives, 
\begin{equation}\label{en30}
I_{k+1} = \frac{m_1! m_2! ... m_{k+1}!}{(N+k)!} 
\end{equation}
For non-integral values of $m_i$ the Eq. \ref{en30} must be remain
expressed in terms of Gamma functions.  
Note that the above expression (\ref{en30}) is for $n=k+1$, and
usually we will  
evaluate it with $n=k$, for which case $I_k = m_1!m_2!
... m_k!/(N+k-1)!$.

To obtain the $q$th moment of $p_i$ one simply needs to substitute $(m_i+q)$ 
for $m_i$ in $I_{k}$, and calculate the ratio of
$I_{k}(m_i+q)/I_{k}(m_i)$, whose meaning is hopefully clear. 
For example, $\langle p_i  \rangle$ is given by, 
\begin{equation}\label{av1}
\langle p_i \rangle = \frac{m_1! m_2! ... (m_i+1)!
  ... m_{k}!}{(N+k)!}  
\frac{(N+k-1)!}{m_1!m_2! ... m_{k}!} 
= \frac{m_i+1}{N+k}  
\end{equation}
where the notation  $\langle p_i \rangle$ is used to denote the
moment of $p_i$ when there are $k$ ``bins''. 
Similarly, 
\begin{equation}
\langle p_i^2 \rangle = \frac{m_1! m_2! ... (m_i+2)!
  ... m_{k}!}{(N+k+1)!}  
\frac{(N+k-1)!}{m_1!m_2! ... m_{k}!} 
= \frac{(m_i+2)(m_i+1)}{(N+k+1)(N+k)}  
\end{equation}
Giving the standard deviation as,
\begin{equation}
\langle p_i^2  \rangle - \langle p_i \rangle^2 
= \frac{(m_i +1)(N+k-m_i-1)}{(N+k)^2 (N+k+1) }
\end{equation}
These results are in agreement with those of Friedman. 
Higher order moments are also easily calculated. 
The difference of the skewness from zero for example, can give an
indication of the extent to which noise in the data should be regarded
as non-Gaussian.  
Note that because $\sum_{i=1}^n p_i = 1$, then, 
\begin{equation}\label{sumtest}
\begin{array}{ll}
1 &= \int_D dp_1 ... dp_{k-1} 
\left( \sum_{i=1}^k p_i \right) 
P(p_1, ... , p_k | m_1, ... , m_k)
\\
& = \sum_{i=1}^k \langle p_i \rangle 
\end{array}
\end{equation}
where $D$ is used as shorthand to indicate that the integral should be
over the correct domain of integration subject to the constraint of
$\sum_{i=1}^n p_i =1$.  
Eq. \ref{sumtest} is correctly satisfied by Eq. \ref{av1}.

\section{Remarks}

There are a variety of distributions in which 
the $\{p_i\}$ only appear in 
a factor of
$\Pi_{i=1}^n p_i^{m_i}$, and the results here apply to those
cases also. 
More generally the probability distribution or its prior could involve
any function of $\{ p_i \}$. 
For example, we might want to introduce a suitable prior into the
problem so as to bias against "outliers", or towards
a particular set of $\{p_i\}$.   
In these more general cases the change of variables
to $n$-dimensional spherical polars will still allow a comparatively
straightforward numerical integral. 
A numerical integral over the $\{ p_i \}$ subject to $0\leq p_i \leq
1$ and $\sum_{i=1}^n p_i=1$, without the change of variables to
spherical polars, is not so easy.
For  
%%Following the change of variables, for 
some combinations of priors and 
probability distributions  the integral will
remain factorisable after the change of variables. 
This might continue to be useful for other analytical 
calculations.

\end{document}